\documentclass{article}
\usepackage{amssymb}
\usepackage{amscd}
\usepackage{amstext}
\usepackage{amsthm}
\usepackage{amsmath}
\usepackage{latexsym}
\title{}
\begin{document}

\theoremstyle{plain}
\newtheorem{theorem}{Theorem}[section]
\newtheorem{corollary}[theorem]{Corollary}
\newtheorem{proposition}[theorem]{Proposition}
\newtheorem{lemma}[theorem]{Lemma}
\newtheorem{definition}[theorem]{Definition}

\newtheorem{remark}{Remark}[section]
\newtheorem{example}{Example}[section]
\newtheorem{exercise}{Exercise}[section]

\title{Eisenstein-Dumas criterion and the action of
 $2\times 2$ nonsingular triangular matrices on polynomials in one variable}
\author{Martin Jur\' a\v s
}

\date{} \maketitle{}


\bigskip

\noindent {\bf Abstract} Let $K$ be a valued field (in general $K$ is not heselian) with valuation $v$ and $A(x)\in K[x]$ be a polynomial of degree $n$.
We find necessary and sufficient conditions for the existence
of the elements $s,t,u\in K$, $s\neq 0\neq u$, such that  at least one of the polynomials
$u^nA(\frac{sx+t}{u})$, $(tx+u)^nA(\frac{sx}{tx+u})$, $(ux)^nA(\frac{tx+s}{ux})$ or
$(ux+t)^nA(\frac{s}{ux+t})$ is an Eisenstein-Dumas polynomial at $v$, provided that
the characteristic of  the residue field of $v$  does not divide $n$. Furthermore, we show that if the orbit $A(x)GL(2,K)$
contains an Eisenstein-Dumas polynomial at $v$, then an Eisenstein-Dumas polynomial at $v$ can be found
in a certain one-parameter subset of this orbit.

\bigskip


%
%

\section{Introduction and results}

In the past 150 years, original criteria of Sch\" onemann \cite{schonemann46} and Eisenstein \cite{eisen50}
were generalized in many different ways. One far reaching generalization comes in the form of a geometric criterion of Dumas \cite{dumas06}:

\medskip

{\it  Newton polygon of a product of polynomials is composed from the sides of polygons of its factors which are ordered
with respect to increasing slopes.}

\medskip

In case when the Newton polygon of a polynomial consists of only one line with no interior lattice points, Dumas' criterion implies
irreducibility. This special case is called Eisenstein-Dumas criterion:

\smallskip
{\it
Let $F$ be a field of fractions of a unique factorization domain $R$. Let  $A(x)=\sum_{i=0}^n a_ix^i \in F[x]$ be a polynomial of degree $n$
and let $p\in R$ be a prime. Denote by $v_p$ the $p$-adic valuation on $F$. If

\medskip

\noindent (D0) \ $a_0a_n\neq 0,$

\medskip

\noindent (D1)  \ $\gcd(v_p(a_0)-v_p(a_n),n)=1$ \ and

\medskip

\noindent (D2) \ $v_p(a_i) \geqslant \frac{n-i}{n} v_p(a_0)+\frac{i}{n} v_p(a_n)
 \quad \text{for} \quad 0\leqslant i\leqslant n,$

\medskip

\noindent then $A(x)$ is irreducible in $F[x].$
}

\medskip

Below, we present a generalization of this criterion to valued fields (see \cite{brown08} or \cite{bishnoi10}) .
For convenience, we first recall some basic facts  about valuations
(see, for instance, Jacobson \cite{jacobson09}  paragraph 9.6).

Let $\Gamma=(\Gamma,+, \leqslant)$ be a {\it linearly ordered abelian group}.
Let $\infty\not\in \Gamma$ and extend the group operation to
$\Gamma \cup\{\infty\}$ by $g+\infty=\infty+g=\infty+\infty=\infty$  and set $g< \infty$, for all $g\in\Gamma$.
Let $K$ be a field and let $K^{\times}$ denotes the multiplicative group of $K$. A function
$v: K\rightarrow\Gamma \cup\{\infty\}$ is referred to as {\it Krull valuation of $K$}, or {\it valuation of $K$},
if it satisfies the following conditions:

\medskip

\noindent (i) \ \ $v(a)=\infty$ iff $a=0$, for all $a\in K$;

\noindent (ii) \ $v(ab)=v(a)+v(b)$, for all $a,b\in K$;

\noindent (iii) $v(a+b)\geqslant\min\{v(a),v(b)\}$, for all $a,b\in K$.

\medskip

\noindent A field endowed with a valuation is called a {\it valued field}.

The ring  $R_v=\{ a\in K: v(a) \geqslant 0\}$ is called the {\it valuation ring of $v$}.
The ideal $M_v=\{ a\in K: v(a) > 0\}$ is called the {\it prime ideal of $v$} or the {\it maximal ideal of $v$}.
The field $K_v=R_v/M_v$ is called the  {\it residue field of $v$}.
Any valuation $v$  satisfies these basic properties:  $v(1)=v(-1)=0$,
 $v(a^{-1})=-v(a)$  for all $a\in K^{\times}$, and $v(a+b)=\min\{v(a),v(b)\}$ provided $v(a)\neq v(b)$.

\medskip

{\it Throughout this paper, $K$ denotes a valued field with value group  $\Gamma$ and  a non-trivial valuation $v:K\rightarrow \Gamma\cup\{\infty\}$.
$K$ is not assumed to be henselian, unless specifically stated.}

\medskip

Theorem \ref{T.1.1} is a special
case of a result by  Brown (cf. \cite{brown08}, Lemma 4).

\begin{theorem}\label{T.1.1} [Eisenstein-Dumas criterion]
 Let
$A(x)=\sum_{i=0}^n a_ix^i\in K[x]$ be a polynomial of degree $n$. If
\medskip

\noindent (D0) \ $a_na_0\neq0$;
\medskip

\noindent (D1) \  $v(a_0)-v(a_n)\not\in k\Gamma=\{kg : g\in \Gamma\}$ for any integer $k>1$ that divides $n$;

\medskip

\noindent (D2) \ $nv(a_i)\geqslant (n-i) v(a_0)+i v(a_n)$ \quad for \quad $0\leqslant i\leqslant n$;

\medskip

\noindent then $A(x)$ is irreducible in $K[x].$
\end{theorem}

\bigskip

A polynomial satisfying  conditions (D0), (D1) and (D2) of Theorem \ref{T.1.1} is referred to as  {\it Eisenstein--Dumas  polynomial (at $v$).}

Classical invariant theory studies intrinsic properties of binary homogeneous forms
 $A(x,y)=\sum_{i=0}^n a_ix^iy^{n-i}$, $a_i\in K$, with regards to the right action of $GL(2,K)$:
$$
A(x,y)g=\sum_{i=0}^n a_i(ax+by)^i(cx+dy)^{n-i}
\quad \text {for} \ \
 g=\left[
       \begin{array}{cc}
         a & b \\
         c & d \\
       \end{array}
     \right]
\in GL(2,K).
$$
Each binary form of degree $n$, corresponds to the (inhomogeneous) polynomial $A(x)=A(x,1)$.
Conversely, given a polynomial $A(x)$, we can recover its homogeneous form by the rule
$A(x,y)=y^nA(\frac{x}{y})$, provided we a priori  specify the degree $n$. The degree of a non-trivial binary form is the degree of any of its monomials.
{\it The degree of a polynomial $A(x)$ is equal to the degree of the homogeneous form corresponding  to $A(x)$} (cf., for instance, \cite{olver99}, Chapter 2).

The action of $GL(2,K)$
on $K[x]$ is given by the linear fraction transformation:
$$
A(x)g=(cx+d)^nA\left(\frac{ax+b}{cx+d}\right) \quad\text{for} \ \ A(x)\in K[x], \ \
g=\left[
       \begin{array}{cc}
         a & b \\
         c & d \\
       \end{array}
     \right]
\in GL(2,K),
$$
where $n$ is the degree  of $A(x)$.

The definition of the degree above guaranties that {\it all polynomials in the orbit $A(x)GL(2,K)$ have the same degree.}
For instance, consider the $3^\text{rd}$-degree polynomial $A(x)=x^3+x^2-2$,
which corresponds to
the form $A(x,y)=x^3+x^2y-2y^3$.
Then,
$
B(x)=A(x)\left[
\begin{array}{cc}
1 & 0 \\
 1 & 1 \\
\end{array}
\right]= -5x^2-6x-2
$
is a polynomial of degree 3, corresponding to the binary form $B(x,y)=-5x^2y-6xy^2-2y^3$.

\medskip

{\it Given a polynomial $A(x)\in K[x]$ of degree $n$, we are interested in finding criteria which ensure
that the orbit $A(x)GL(2,K)$ contains an Eisenstein-Dumas polynomial.}

\medskip

In this paper we obtained some partial results.

Note that all irreducible polynomials cannot be found  by application Eisenstein-Dumas criterion and the action of $GL(2,K)$.
For example,  $A(x)=x^4-14x^2+9$  is irreducible in  $\mathbb Q[x]$, but the orbit $A(x)GL(2,\mathbb Q)$ contains no
Eisenstein-Dumas polynomial at any  valuation $v$ of $\mathbb Q$ (cf. \cite{juras06}, page 229).

In \cite{juras06}, Theorem 2.3, we proved:

\medskip

{\it
Let $F$ be a field of fractions of a unique factorization domain $R$. Let  $A(x)=\sum_{i=0}^n a_ix^i \in F[x]$ be a polynomial of degree $n$
and let $p\in R$ be a prime. Assume that the characteristic of $R$ does not divide $n$  and $p$ does not divide
$n1_R$. If there are $s,t\in F,$
$s\neq 0$ such that $A(sx+t)$ is an Eisenstein-Dumas polynomial at $v_p$, then $A(x-\frac{a_{n-1}}{na_n})$ is an Eisenstein-Dumas
polynomial at the $p$-addic valuation $v_p$.
}

\medskip

A straightforward application of this last theorem  yields the irreducibility,
in $\mathbb Q[x]$, of the $p$-th cyclotomic polynomial $\Phi_p(x)=\sum_{i=0}^{p-1} x^i$ (for an odd prime $p$);
by simple divisibility arguments, one can show that $\Phi_p(x-\frac{1}{p-1})$ is an Eisenstein
polynomial at  $v_p$.
Notice that the transformation suggested by the theorem, $x\mapsto x-\frac{1}{p-1}$, is not the usual one
$x\mapsto x+1$. But, $-\frac{1}{p-1}=1$  in $GF(p)$, and the two transformations are connected
 by  means of Lemma \ref{L.2.2}, since $\Phi_p(x-\frac{1}{p-1})=\Phi_p(x+1-\frac{p}{p-1})$ and $v(-\frac{p}{p-1})=1>\frac{1}{p-1}$.

 Bishnoi and Khanduja  generalized \cite{juras06}, Theorem 2.3,   to henselian valued fields (\cite{bishnoi10}, Theorem 1.2):

\medskip

{\it Let  $g(x)=\sum_{i=1}^ea_ix^i$ be a monic polynomial with coefficients in a henselian field $(K,v)$.
Suppose that the characteristic of the residue field of $v$ does not divide $e.$ If there exists an element
$b\in K$ such that $g(x+b)$ is an Eisenstein-Dumas polynomial with respect to $v$, then so is $g(x-\frac{a_{e-1}}{e})$.}

\medskip

The next result generalizes the theorem above.

\begin{theorem}\label{T.1.0}
Let $A(x)=\sum_{i=0}^n a_ix^i\in K[x]$ be a polynomial of degree $n$.
Assume that the characteristic of the residue field of $v$  does not divide $n$.
If, for some element
$g=\left[
     \begin{array}{cc}
       s & t \\
       u & v \\
     \end{array}
   \right]
\in GL(2,K)$ with $stuv=0$, $A(x)g$ is an Eisenstein-Dumas polynomial, then
at least one of the two polynomials
$$
U(A(x))=A(x)\left[
        \begin{array}{cc}
          1 & -\frac{a_{n-1}}{na_n} \\
          0 & 1 \\
        \end{array}
\right]
\quad \text{or} \quad
L(A(x))=A(x)\left[
        \begin{array}{cc}
          1 & 0 \\
          -\frac{a_{1}}{na_0} & 1 \\
        \end{array}
      \right]
$$
is an Eisenstein-Dumas polynomial.
\end{theorem}

In the proof of \cite{bishnoi10}, Theorem 1.2, Bishnoi and Khanduja  utilized ramification theory of algebraic field extensions
and used the fact that the field $K$ is henselian.
Theorem \ref{T.1.0} applies to {\it any} valued field $K$; our proof is a modification of the author's proof of \cite{juras06}, Theorem 2.3.

Let $A(x)\in \mathbb Q[x]$ be a polynomial of degree $n$.
Theorem \ref{T.1.0}  implies that there are  at most finitely many prime numbers $p\in\mathbb Z$, for which there is
$
g=\left[
\begin{array}{cc}
   s & t \\
       u & v \\
        \end{array}
\right]\in GL(2,\mathbb Q)
$,
$stuv=0$, such that $A(x)g$ is an Eisenstein-Dumas polynomial at some $p$-adic valuation $v_p$. Indeed, for such a $p$, if $b_0$ and $b_n$ denote
the constant coefficient and the leading coefficient of $U(A(x))$, and $c_0$ and $c_n$ denote
the constant coefficient and the leading coefficient of $L(A(x))$, then, due to (D1),
 $p$ satisfies at least one of the three conditions below:
(i) $\gcd(v_p(\frac{b_0}{b_n}), n)=1$  or (ii) $\gcd(v_p(\frac{c_0}{c_n}), n)=1$  or (iii) $p$ divides $n$.

\begin{theorem} \label{T.1.8}
Let $A(x)=\sum_{i=0}^n a_ix^i\in K[x]$ be a polynomial of degree $n$.
If one of the two polynomials
$U(A(x))$ or $L(A(x))$
is an Eisenstein-Dumas polynomial, then $A(x)$ is irreducible.
\end{theorem}

Theorem \ref{T.1.8} follows immediately from Theorem \ref{T.1.1} and the fact that the incredibility of $A(x)$ is invariant under
the action of $GL(2,K)$.
Due to Theorem \ref{T.1.0}, one may expect Theorem \ref{T.1.8} to be more useful in establishing irreducibility of polynomials than Theorem \ref{T.1.1}.
There are instances though, when Theorem  \ref{T.1.1} applies and
Theorem \ref{T.1.8} fails. This may occur when  the characteristic
of  the residue field of $v$  divides the degree of $A(x)$. For example, consider  $A(x)=x^2+4x+8\in \mathbb  Q[x]$, which is an Eisenstein-Dumas
polynomial at the 2-adic valuation  $v_2$, but neither $U(A(x))=x^2+4$ nor $L(A(x))=\frac{1}{2} x^2+8$ is an
Eisenstein-Dumas polynomial.\footnote{Condition (D1) of Theorem \ref{T.1.1} fails in both cases.}

\begin{proposition} \label{T.1.9}
Let $A(x)=\sum_{i=0}^n a_ix^i\in K[x]$ be a polynomial of degree $n$  and let
$B(x)=\sum_{i=0}^n a_{n-i}x^i$. Then
$$
U(B(x))=L(A(x))\left[
                 \begin{array}{cc}
                   0 & 1 \\
                   1 & 0 \\
                 \end{array}
               \right].
$$
Furthermore, $L(A(x))$ is an Eisenstein-Dumas polynomial iff \, $U(B(x))$
is an Eisenstein-Dumas polynomial.
\end{proposition}

\begin{theorem}\label{T.1.7}
Let $A(x)=\sum_{i=0}^n a_ix^i\in K[x]$ be a polynomial of degree $n$.
Assume that the characteristic of  the residue field of $v$  does not divide $n$.
If the orbit  $A(x)GL(2,K)$ contains an Eisenstein-Dumas polynomial,  then one of the two polynomials $U(A(x))$ or
$L(A(x))$ is an Eisenstein-Dumas polynomial or there is $t\in K$ such that $A'(t)\neq 0$ and the polynomial
$$
A(x)\left[
      \begin{array}{cc}
        t & \phi(t) \\
        1 & 1 \\
      \end{array}
    \right]
=(x+1)^n\,A\left(\frac{tx+\phi(t)}{x+1}\right),
\quad \text{where} \quad \phi(t)=t-n\frac{A(t)}{A'(t)},
$$
is an  Eisenstein-Dumas polynomial.
($A'$ is the derivative of $A$.)
\end{theorem}

This last result generalizes \cite{juras06}, Theorem 3.2.

%
%

\section{Proofs}

Theorem \ref{T.1.1} not only generalizes numerous versions of Sch\" onemann-Eisenstein criterion, but also exhibits
 symmetries of the irreducibility conditions that are absent under the conditions of the primer criteria of Sch\" onemann
 \cite{schonemann46} and Eisenstein \cite{eisen50}.  There are
two obvious families of symmetries generated by group actions of $K^\times$  and a discrete symmetry.
The proof of the following lemma is straightforward.

\begin{lemma} \label{L.1.2}
Let $A(x)\in K[x]$ and $t\in K^\times.$
Then

\medskip

\noindent (i) \ \ $tA(x)$ is an Eisenstein--Dumas  polynomial  iff $A(x)$ is an Eisenstein--Dumas  polynomial.

\noindent (ii) \ $A(tx)=A(x)\left[
                              \begin{array}{cc}
                                t & 0 \\
                                0 & 1 \\
                              \end{array}
                            \right]
$
is an Eisenstein--Dumas  polynomial  iff $A(x)$ is an Eisenstein--Dumas  polynomial.

\noindent (iii)  $x^nA(\frac{1}{x})\left[
                                    \begin{array}{cc}
                                      0 & 1 \\
                                      1 & 0 \\
                                    \end{array}
                                  \right]
$
is an Eisenstein--Dumas  polynomial iff $A(x)$ is an Eisenstein--Dumas  polynomial.

\end{lemma}

\begin{lemma}\label{L.2.1}
Let
$A(x)=\sum_{i=0}^n a_ix^i\in K[x]$ be a polynomial of degree $n$.
Then the conditions (D0),(D1),(D2) in Theorem \ref{T.1.1} are equivalent to (D0),(D1), (D2'), where

\medskip

\noindent (D2') \ $nv(a_i)> (n-i) v(a_0)+i v(a_n)$ for  $1\leqslant i\leqslant n-1$.
\end{lemma}
\proof%
(D2) is always trivially satisfied for $i=0$ and $i=n$ and so (D0),(D1),(D2') implies (D0),(D1),(D2). Let (D0),(D1), (D2) hold and assume, by contraposition,
that $nv(a_i)=(n-i) v(a_0)+i v(a_n)$, for some $i$,  $1\leqslant i\leqslant n-1$. Then $ig=nh$ where $g=v(a_0)-v(a_n)$ and $h=v(a_0)-v(a_i)$. Let $c=\gcd(i,n)$.
 Obviously, $c<n$. There are $a,b\in \mathbb Z$ such that $ai+bn=c$. Multiplying the last equation by  $\frac{g}{c}$ and using $ig=nh$ we arrive
 at $\frac{n}{c}(ah+bg)=g$. $\frac{n}{c}$ is an integer greater than 1 which divides $n$ and $ah+bg\in \Gamma$.
The last equation contradicts (D1).
\hfill $\rule{0.5em}{0.5em}$

\begin{lemma} \label{L.2.2}
Let $A(x)=\sum_{i=0}^n a_ix^i\in K[x]$ be an Eisenstein-Dumas polynomial of degree $n$.
Assume that $t\in K$, $nv(t) > v(a_0)-v(a_n)$.
Then
$$
A(x)\left[
        \begin{array}{cc}
          1 & t \\
          0 & 1 \\
        \end{array}
      \right]=A(x+t)
$$
is an  Eisenstein-Dumas polynomial.
\end{lemma}
\proof
$
\bar A(x)=\sum_{i=0}^n \bar a_ix^i=A(x+t),
$
where
$$
\bar a_i=\sum_{k=i}^n \left(
                        \begin{array}{c}
                          k \\
                          i \\
                        \end{array}
                      \right) a_kt^{k-i}  \qquad \text{for all} \quad 0\leqslant i \leqslant n.
$$
In particular, $\bar a_n=a_n\neq 0$, and so $v(\bar a_n)=v(a_n)$, and $\bar a_0=A(t)$. If $\bar a_0=0$,
then $t$ is a root of $A(x)$, a contradiction. Thus $\bar a_0\neq 0$ and so (D0) is satisfied for $\bar A(x)$.

For $0<k\leqslant n$, we have
$$
n v(a_kt^k)=n[ v(a_k)+kv(t)]\geqslant (n-k)v(a_0) +k v(a_n) +kn v(t)
$$
$$
> (n-k)v(a_0) +k v(a_n) + k(v(a_0)-v(a_n))= n v(a_0).
$$
Hence, $v(a_kt^k)>v(a_0)$ for $0<k\leqslant n$.
Since, $\bar a_0=a_0+\sum_{k=1}^{n} a_kt^k$, then  $v(\bar a_0)=v(a_0)$
and so (D1) holds for $\bar A(x).$
\medskip

We now show that (D2) is satisfied for $\bar A(x)$. For $i\leqslant k,$ we have
$$
n v(\left( \left(
               \begin{array}{c}
                 k \\
                 i \\
               \end{array}
             \right) a_kt^{k-i}\right)
             \geqslant n[v(a_k)+(k-i)v(t)]
\geqslant (n-k)v(a_0) + k v(a_n)
$$
$$
+(k-i)(v(a_0)-v(a_n))
= (n-i) v(a_0) +i v(a_n) =(n-i) v(\bar a_0) +i v(\bar a_n)
$$
Thus, for $0\leqslant i\leqslant n$,
$$
n v(\bar a_i)=n v\left(\sum_{k=i}^n \left(
                        \begin{array}{c}
                         k \\
                         i \\
                        \end{array}
                      \right) a_kt^{k-i}\right)\geqslant (n-i) v(\bar a_0) +i  v(\bar a_n).
$$
\hfill $\rule{0.5em}{0.5em}$

\bigskip

To prove the following lemma, we need   $v(n1)=0$. Bishnoi and Khanduja (see \cite{bishnoi10}, Theorem 1.2)
found an elegant equivalent condition, namely that the characteristic of the residue field of $v$  does not divide $n$.

\begin{lemma} \label{L.2.3}
Let $A(x)=\sum_{i=0}^n a_ix^i\in K[x]$ be a polynomial of degree $n$.
Assume that the characteristic of the residue field of $v$  does not divide $n$.
If $A(x)$ is an Eisenstein-Dumas polynomial, then $nv(\frac{a_{n-1}}{na_n}) > v(a_0)-v(a_n)$.
\end{lemma}
\proof
First we use the fact that the characteristic of the residue field of $v$, $K_v$,  does not divide $n$
to show that  $v(n1)=0$.

Assume, by contraposition, that $v(n1)\neq 0$. Since $v(1)=0$, it follows $v(n1)> 0$, and so
$n1\in M_v$. Hence, $n\varphi(1)=\varphi(n1)=0$, where $\varphi: R_v\rightarrow R_v/M_v$ is the canonical homomorphism.
Thus,  the characteristic of $K_v$ divides $n$,
which yields a contradiction.

Since $n1\neq 0$, then $na_n=(n1)a_n\neq 0$.

Furthermore, assume that $A(x)$ is an Eisenstein-Dumas polynomial.
Using (D2') in Lemma \ref{L.2.1} with $i=n-1$, we obtain
$$
n v(\frac{a_{n-1}}{n a_n})= n [v(a_{n-1}) - v(a_n)-v(n1)]
$$
$$
>  v(a_0)+(n-1)v(a_n) -n v(a_n)
=v(a_0) - v(a_n).
$$
\hfill $\rule{0.5em}{0.5em}$

\bigskip

\noindent{\bf Proof of Theorem \ref{T.1.0}}

Theorem \ref{T.1.0} will be established by proving Theorems \ref{T.1.3} -- \ref{T.1.6}. Each of these theorems
will deal with a partial case of Theorem \ref{T.1.0}.

By $U(2,K)$, resp., $L(2,K)$ we denote the multiplicative group of all
upper, resp., lower triangular matrices in $GL(2,K).$

\begin{theorem}\label{T.1.3}
Let $A(x)=\sum_{i=0}^n a_ix^i\in K[x]$ be a polynomial of degree $n$.
Assume that the characteristic of the residue field of $v$  does not divide $n$.
If the orbit  $A(x)U(2,K)$ contains an Eisenstein-Dumas polynomial, then
\begin{equation}
U(A(x))=A(x)\left[
        \begin{array}{cc}
          1 & -\frac{a_{n-1}}{na_n} \\
          0 & 1 \\
        \end{array}
\right]=A\left(x-\frac{a_{n-1}}{na_n}\right) \label{E1}
\end{equation}
is an Eisenstein-Dumas polynomial.
\end{theorem}
\proof
Let
$
g=\left[
         \begin{array}{cc}
           s & t \\
           0 & u \\
         \end{array}
       \right]\in U(2,K)
$
and assume that
$$
\bar A(x)=\sum_{i=0}^n \bar a_ix^i=A(x)g=u^nA(\frac{sx+t}{u})
$$
is an Eisenstein-Dumas polynomial.

By Lemma \ref{L.2.3},
$
n v(\frac{\bar a_{n-1}}{n\bar a_n})> v(\bar a_0) - v(\bar a_n)
$
 and by Lemma \ref{L.2.2},
$$
B(x)=\bar A(x)\left[
        \begin{array}{cc}
          1 & -\frac{\bar a_{n-1}}{n\bar a_n} \\
          0 & 1 \\
        \end{array}
      \right]
$$
is an Eisenstein-Dumas polynomial. By Lemma \ref{L.1.2}(i,ii),
$(\frac{1}{u})^nB(\frac{u}{s}x)$ is an Eisenstein-Dumas polynomial.
Using  $0\neq\bar a_n=a_ns^n$ and $\bar a_{n-1}=s^{n-1}(na_nt+a_{n-1}u)$,  we obtain
$$
(\frac{1}{u})^nB(\frac{u}{s}x)=B(x)\left[
                        \begin{array}{cc}
                          \frac{1}{s} & 0 \\
                          0 & \frac{1}{u} \\
                        \end{array}
                      \right]
=\bar A(x)\left[
        \begin{array}{cc}
          1 & -\frac{\bar a_{n-1}}{n\bar a_n}  \\
          0 & 1 \\
        \end{array}
      \right]
      \left[
                        \begin{array}{cc}
                          \frac{1}{s} & 0 \\
                          0 & \frac{1}{u} \\
                        \end{array}
                      \right]
$$
$$
= A(x)\left[
        \begin{array}{cc}
          s & t \\
          0 & u \\
        \end{array}
      \right]
      \left[
        \begin{array}{cc}
          1 & -\frac{t}{s} -\frac{a_{n-1}}{na_n}\frac{u}{s} \\
          0 & 1 \\
        \end{array}
      \right]
\left[
                        \begin{array}{cc}
                          \frac{1}{s} & 0 \\
                          0 & \frac{1}{u} \\
                        \end{array}
                      \right]
=A(x)       \left[
        \begin{array}{cc}
          1 & -\frac{a_{n-1}}{n a_n}  \\
          0 & 1 \\
        \end{array}
      \right].
$$
\hfill $\rule{0.5em}{0.5em}$

%
%
%

\begin{theorem}\label{T.1.4}
Let $A(x)=\sum_{i=0}^n a_ix^i\in K[x]$ be a polynomial of degree $n$.
Assume that the characteristic of  the residue field of $v$  does not divide $n$.
If the orbit  $A(x)L(2,K)$ contains an Eisenstein-Dumas polynomial,
then
\begin{equation}
L(A(x))=A(x)\left[
        \begin{array}{cc}
          1 & 0 \\
          -\frac{a_{1}}{na_0} & 1 \\
        \end{array}
      \right]=(1-\frac{a_1}{na_0}x)^nA\left(\frac{x}{1-\frac{a_1}{na_0}x}\right) \label{E3}
\end{equation}
is an Eisenstein-Dumas polynomial.
\end{theorem}
\proof
Let
$
g=\left[
         \begin{array}{cc}
           s & 0 \\
           t & u \\
         \end{array}
       \right]\in L(2,K)
$
Assume that
$
\bar A(x)=A(x)g
$
is an Eisenstein-Dumas polynomial.
By Lemma \ref{L.1.2}(iii), $x^n\bar A(\frac{1}{x})$ is also an Eisenstein-Dumas polynomial.
We have
$$
x^n\bar A(\frac{1}{x})=\bar A(x)\left[
                \begin{array}{cc}
                  0 & 1 \\
                  1 & 0 \\
                \end{array}
              \right]
=A(x)\left[
         \begin{array}{cc}
           s & 0 \\
           t & u \\
         \end{array}
       \right]\left[
                \begin{array}{cc}
                  0 & 1 \\
                  1 & 0 \\
                \end{array}
              \right]
$$
$$
=A(x)\left[
                \begin{array}{cc}
                  0 & 1 \\
                  1 & 0 \\
                \end{array}
              \right]\left[
\begin{array}{cc}
           u & t \\
           0 & s \\
         \end{array}
       \right]\left[
                \begin{array}{cc}
                  0 & 1 \\
                  1 & 0 \\
                \end{array}
              \right]\left[
                \begin{array}{cc}
                  0 & 1 \\
                  1 & 0 \\
                \end{array}
              \right]
=B(x)\left[
\begin{array}{cc}
           u & t \\
           0 & s \\
         \end{array}
       \right]
$$
where
$
B(x)=A(x)\left[
                \begin{array}{cc}
                  0 & 1 \\
                  1 & 0 \\
                \end{array}
              \right]=\sum_{i=0}^n a_{n-i}x^i.
$
Since
$
\left[
\begin{array}{cc}
           u & t \\
           0 & s \\
         \end{array}
       \right]\in U(2,K),
$
then by Theorem \ref{T.1.3},
$
B(x)\left[
      \begin{array}{cc}
        1 & -\frac{a_1}{na_0} \\
        0 & 1 \\
      \end{array}
    \right]
$
is an Eisenstein-Dumas polynomial.
$$
B(x)\left[
      \begin{array}{cc}
        1 & -\frac{a_1}{na_0} \\
        0 & 1 \\
      \end{array}
    \right]
=A(x)\left[
                \begin{array}{cc}
                  0 & 1 \\
                  1 & 0 \\
                \end{array}
              \right]
\left[
      \begin{array}{cc}
        1 & -\frac{a_1}{na_0} \\
        0 & 1 \\
      \end{array}
    \right]
=A(x)\left[
       \begin{array}{cc}
         0 & 1 \\
         1 & -\frac{a_1}{na_0} \\
       \end{array}
     \right]
$$
By Lemma \ref{L.1.2}(iii)
$$
A(x)\left[
       \begin{array}{cc}
         0 & 1 \\
         1 & -\frac{a_1}{na_0} \\
       \end{array}
     \right]\left[
              \begin{array}{cc}
                0 & 1 \\
                1 & 0 \\
              \end{array}
            \right]
            =A(x)\left[
                   \begin{array}{cc}
                     1 & 0 \\
                     -\frac{a_1}{na_0} & 1 \\
                   \end{array}
                 \right]=L(A(x))
$$
is an Eisenstein-Dumas polynomial.
\hfill $\rule{0.5em}{0.5em}$

\begin{theorem}\label{T.1.5}
Let $A(x)=\sum_{i=0}^n a_ix^i\in K[x]$ be a polynomial of degree $n$.
Assume that the characteristic of  the residue field of $v$  does not divide $n$.
If the orbit,
$
A(x)
U(2,K)
   \left[
     \begin{array}{cc}
       0 & 1 \\
       1 & 0 \\
     \end{array}
   \right]
$
contains an Eisenstein-Dumas polynomial,
then the polynomial $U(A(x))$
is an Eisenstein-Dumas polynomial.
\end{theorem}
\proof
Let
$g=\left[
     \begin{array}{cc}
       t & s \\
       u & 0 \\
     \end{array}
   \right]\in U(2,K)
   \left[
     \begin{array}{cc}
       0 & 1 \\
       1 & 0 \\
     \end{array}
   \right]
$
be such that
$
B(x)=A(x)g=(ux)^nA\left(\frac{tx+s}{ux}\right)
$
is an Eisenstein-Dumas polynomial.

By Lemma \ref{L.1.2}(iii),
$x^nB(\frac{1}{x})$ is an Eisenstein-Dumas polynomial. We have
$$
x^nB\left(\frac{1}{x}\right)=B(x)\left[
       \begin{array}{cc}
         0 & 1 \\
         1 & 0 \\
       \end{array}
     \right]
=A(x)\left[
       \begin{array}{cc}
         t & s \\
         u & 0 \\
       \end{array}
     \right]
     \left[
       \begin{array}{cc}
         0 & 1 \\
         1 & 0 \\
       \end{array}
     \right]
=A(x)\left[
       \begin{array}{cc}
         s & t \\
         0 & u \\
       \end{array}
     \right]
$$
Since
$
\left[
       \begin{array}{cc}
         s & t \\
         0 & u \\
       \end{array}
     \right]
\in U(2,K)$, the statement  follows from Theorem \ref{T.1.3}.
\hfill $\rule{0.5em}{0.5em}$

\begin{theorem}\label{T.1.6}
Let $A(x)=\sum_{i=0}^n a_ix^i\in K[x]$ be a polynomial of degree $n$.
Assume that the characteristic of  the residue field of $v$  does not divide $n$.
If the orbit
$
A(x)L(2,K)\left[
                       \begin{array}{cc}
                         0 & 1 \\
                         1 & 0 \\
                       \end{array}
                     \right]
$
contains an Eisenstein-Dumas polynomial,
then the polynomial $L(A(x))$
is an Eisenstein-Dumas polynomial.
\end{theorem}
\proof%
Let
$
g=\left[
      \begin{array}{cc}
        0 & s \\
        u & t \\
      \end{array}
    \right]\in L(2,K)\left[
                       \begin{array}{cc}
                         0 & 1 \\
                         1 & 0 \\
                       \end{array}
                     \right]
$
be such that $B(x)=A(x)g=(ux+t)^nA(\frac{s}{ux+t})$ is an Eisenstein-Dumas polynomial.
 By Lemma \ref{L.1.2} (iii),
$x^nB(\frac{1}{x})$ is an Eisenstein-Dumas polynomial of degree $n$. We have
$$
x^nB\left(\frac{1}{x}\right)=B(x)\left[
       \begin{array}{cc}
         0 & 1 \\
         1 & 0 \\
       \end{array}
     \right]
=A(x)\left[
       \begin{array}{cc}
         0 & s \\
         u & t \\
       \end{array}
     \right]
     \left[
       \begin{array}{cc}
         0 & 1 \\
         1 & 0 \\
       \end{array}
     \right]
=A(x)\left[
       \begin{array}{cc}
         s & 0 \\
         t & u \\
       \end{array}
     \right]
$$
Since
$
\left[
       \begin{array}{cc}
         s & 0 \\
         t & u \\
       \end{array}
     \right]
\in L(2,K)$, the statement now follows from Theorem \ref{T.1.4}.
\hfill $\rule{0.5em}{0.5em}$

\bigskip

This completes the proof of Theorem \ref{T.1.0}.

\bigskip

\noindent{\bf Proof of Proposition \ref{T.1.9}}
\proof%
We have
$$
U(B(x))=B(x)\left[
              \begin{array}{cc}
                1 & -\frac{a_1}{na_0} \\
                0 & 1 \\
              \end{array}
            \right]
=A(x)\left[
       \begin{array}{cc}
         0 & 1 \\
         1 & 0 \\
       \end{array}
     \right]\left[
              \begin{array}{cc}
                1 & -\frac{a_1}{na_0} \\
                0 & 1 \\
              \end{array}
            \right]
$$
$$
=A(x)\left[
        \begin{array}{cc}
          0 & 1 \\
          1 & -\frac{a_1}{na_0} \\
        \end{array}
      \right]
=A(x)\left[\begin{array}{cc}
                1 & 0 \\
                -\frac{a_1}{na_0} & 1 \\
              \end{array}
            \right]\left[
                     \begin{array}{cc}
                       0 & 1 \\
                       1 & 0 \\
                     \end{array}
                   \right]
=L(A(x))\left[
                     \begin{array}{cc}
                       0 & 1 \\
                       1 & 0 \\
                     \end{array}
                   \right].
$$
The last sentence now follows from Lemma \ref{L.1.2}(iii).
\hfill $\rule{0.5em}{0.5em}$

\bigskip
\noindent{\bf Proof of Theorem \ref{T.1.7}}
\proof%
Let
$
\left[
  \begin{array}{cc}
    s & t \\
    u & v \\
  \end{array}
\right] \in GL(2,K)
$
and assume  $\bar A(x)=A(x)g$ is an Eisenstein-Dumas polynomial.
If $stuv=0$, then by Theorem \ref{T.1.0},  one of the two polynomials $U(A(x))$ or
$L(A(x))$ is an Eisenstein-Dumas polynomial.

Assume that $stuv\neq 0.$ By Lemma \ref{L.1.2}(i,ii),
$B(x)=\sum_{i=1}^n b_i x^i=(\frac{1}{v})^n\bar A(\frac{v}{u}x)$ is an Eisenstein-Dumas.
$$
B(x)=\bar A(x)\left[
                                \begin{array}{cc}
                                  \frac{1}{u} & 0 \\
                                  0 & \frac{1}{v} \\
                                \end{array}
                              \right]
=A(x)\left[
       \begin{array}{cc}
         s & t \\
         u & v \\
       \end{array}
     \right]
     \left[
                                \begin{array}{cc}
                                  \frac{1}{u} & 0 \\
                                  0 & \frac{1}{v} \\
                                \end{array}
                              \right]
=A(x)\left[
       \begin{array}{cc}
         \frac{s}{u} & \frac{t}{v} \\
         1 & 1 \\
       \end{array}
     \right].
$$
By Lemma \ref{L.2.3},  $nv(\frac{b_{n-1}}{nb_n}) > v(b_0)-v(b_n)$.
By Lemma \ref{L.2.2},
$$
\bar B(x)=B(x)\left[
        \begin{array}{cc}
          1 & -\frac{b_{n-1}}{nb_n} \\
          0 & 1 \\
        \end{array}
      \right]=B(x-\frac{b_{n-1}}{nb_n})
$$
is an  Eisenstein-Dumas polynomial.
Moreover,
$$
b_n=A(e) \qquad \text{and} \qquad b_{n-1}= (f-e)\, A'(e)+n\, A(e),
$$
where $e=\frac{s}{u}$ and $f=\frac{t}{v}$.
$A(e)=b_n\neq 0$ since $B(x)$ is of degree $n$.

If $A'(e)=0$, then
$$
\bar B(x)
=A(x)\left[
       \begin{array}{cc}
         e & f \\
         1 & 1 \\
       \end{array}
     \right]\left[
                \begin{array}{cc}
                  1 & -1 \\
                  0 & 1 \\
                \end{array}
              \right]
=A(x)\left[
       \begin{array}{cc}
         e & f-e \\
         1 & 0 \\
       \end{array}
     \right],
$$
and so by Theorem \ref{T.1.5},  $U(A(x))$ is an Eisenstein-Dumas polynomial.

Assume that $A'(e)\neq 0$. Then, by Lemma \ref{L.1.2}(i,ii),
$$
C(x)=\left(\frac{nA(e)}{(e-f)A'(e)}\right)^n \bar B\left(\frac{(e-f)A'(e)}{nA(e)}x\right)
$$
is an  Eisenstein-Dumas polynomial. We have
$$
C(x)=\bar B(x)\left[
           \begin{array}{cc}
             1 & 0 \\
             0 & \frac{nA(e)}{(e-f)A'(e)} \\
           \end{array}
         \right]
$$
$$
=B(x)\left[
       \begin{array}{cc}
         1 & \frac{(f-e)\, A'(e)+n\, A(e)}{nA(e)} \\
         0 & 1 \\
       \end{array}
     \right]\left[
           \begin{array}{cc}
             1 & 0 \\
             0 & \frac{nA(e)}{(e-f)A'(e)} \\
           \end{array}
         \right]
$$
$$
=A(x)\left[
                \begin{array}{cc}
                  e & f \\
                  1 & 1 \\
                \end{array}
              \right]\left[
       \begin{array}{cc}
         1 & \frac{(f-e)\, A'(e)+n\, A(e)}{nA(e)} \\
         0 & 1 \\
       \end{array}
     \right]\left[
           \begin{array}{cc}
             1 & 0 \\
             0 & \frac{nA(e)}{(e-f)A'(e)} \\
           \end{array}
         \right]
$$
$$
=A(x)\left[
       \begin{array}{cc}
         e & e-n\frac{A(e)}{A'e)} \\
         1 & 1 \\
       \end{array}
     \right]
$$
\hfill $\rule{0.5em}{0.5em}$


\end{document}